# Irreducible Characters for the Symmetric Groups and Kostka Matrices

by

Ronald F. Fox
Smyrna, Georgia
May 12, 2018

### Abstract

In an earlier paper [1] it was shown that the Frobenius compound characters for the symmetric groups are related to the irreducible characters by a linear relation that involves a unitriagular coupling matrix that gives the Frobenius characters in terms of linear combinations of the irreducible characters. It is desirable to invert this relationship since we have formulas for the Frobenius characters and want the values for the irreducible characters. This inversion is straightforward and yields both the irreducible characters but also the coupling matrix that turns out to be the Kostka matrix in the original direction. We show that if the Frobenius monomial identity is applied a modification of it, equation (22), produces a monomial formula that produces the Kostka matrix inverse without involving the characters of either type. However it is a formidable task to execute this procedure for symmetric groups of even modest order. Alternatively the inversion by means of the unitriangular coupling matrix produces the Kostka matrix and the irreducible characters simultaneously and with much less effort than required for the monomial approach. Moreover there is a surprise.

### I. Introduction.

In my first published paper [1], a relationship between the Frobenius compound characters and the irreducible characters for the symmetric groups, $S_n$, was elucidated. Denote a compound character by $\phi^j$ and an irreducible character by $\chi^i$ where $j$ and $i$ are integers that run from 1 to $k$ and $k$ is the number of classes (cycle structures) for the symmetric group, $S_n$, on the integers $1, 2 \ldots n$. For $n = 2\ or\ 3$, $k = 2\ or\ 3$, respectively, but for

larger values of $n$ we find:
$(n, k) = \{(4,5), (5,7), (6,11), (7,15), (8,22), (9,30), (10,42) \ldots \}$.
This numbering scheme associates $k$ with the largest cycle structure, denoted by $(n)$, and 1 with the smallest cycles, the 1-cycles that occur $n$ times, denoted by $(1^n)$. There is a natural linear ordering of the cycle structures that agrees with the ordering of integers from 1 to $k$. {e.g. for $S_4$, $k = 5, 4, 3, 2, 1$ corresponds respectively with the cycle structures $(4)(3,1)(2^2)(2, 1^2)(1^4)$}

The fundamental theorem proved in [1] is:

(1)

$$\phi^k = a_k^k \chi^k$$
$$\phi^{k-1} = a_k^{k-1} \chi^k + a_{k-1}^{k-1} \chi^{k-1}$$
$$\phi^{k-2} = a_k^{k-2} \chi^k + a_{k-1}^{k-2} \chi^{k-1} + a_{k-2}^{k-2} \chi^{k-2}$$
.
.
.
$$\phi^1 = a_k^1 \chi^k + a_{k-1}^1 \chi^{k-1} + a_{k-2}^1 \chi^{k-2} + \cdots + a_1^1 \chi^1$$

The matrix $a_{k-i}^{k-j}$ is labeled by the character names (classes or cycle structures) and the characters are functions of the classes (cycle structures) e.g. $\phi^k(C_l)$. In [1], using the orthonormality of the irreducible characters,

(2)

$$\frac{1}{n!} \sum_\sigma m_\sigma \chi^j(C_\sigma) \chi^i(C_\sigma) = \delta_{ji}$$

it is noted that (1) and (2) imply

(3)

$$a_{k-i}^{k-j} = \frac{1}{n!} \sum_\sigma m_\sigma \phi^{k-j}(C_\sigma) \chi^{k-i}(C_\sigma)$$

in which the sum is over all classes and the number of elements in class $C_\sigma$ is $m_\sigma$. The central feature manifested by eq.(1) is the *triangularity* of the coupling matrix, i.e.

(4)

$$a_{k-i}^{k-j} = 0 \text{ if } k - i < k - j \text{ or } j < i.$$

The triangularity of the coupling matrix permits an easily implemented solution for the irreducible characters in terms of the given Frobenius compound characters.

In 1967 I did not know what the coupling matrices were from a combinatorial point of view, nor did I know how to compute them other than by inverting the triangular relationship exhibited in eq.(1). Such an inversion yielded both the irreducible characters and the coupling matrix elements simultaneously. In 1990, while teaching a group theory class at Georgia Tech, I noticed that Frobenius' construction of monomials not only connected Frobenius compound characters with the irreducible characters but also provided a monomial formula for the $a_{k-i}^{k-j}$ 's, or even better could be converted to directly produce the inverse of this triangular matrix. In 1967 one of my mentors at The Rockefeller University was Gian-Carlo Rota and he had encouraged me to write my paper. In October of 1990 I wrote to him at MIT, to which he had returned around 1969, and asked if the result was new. He said that by then everything was known and could be found in a book by MacDonald [2]. At the time I was into other matters and also found the book a bit difficult to penetrate. Had I persevered I would have seen that my triangular matrices were Kostka matrices [3] and would have found page 111 full of small ($n = 2, 3, 4, 5, 6$) examples of them. In addition a paper appeared as early as 1978 by Foulkes [4] in which he referred to my 1967 paper as containing a rediscovery of the Kostka tables. At that time the internet and its various search engines were not yet making such connections automatic so I remained ignorant of this paper.

A paper by Kerber [5] appeared in 1969 that addressed my request for an interpretation of the triangle matrices. No mention of Kostka was made but two observations were proved by means other than those used in my paper. Kerber noted that for all $j$, $a_{k-j}^{k-j} = 1$, which means the matrices are *unitriangular* and have determinant equal to 1. He also showed that for all $j$, $a_k^{k-j} = 1$ as well. I will verify these observations below using my methods. Kerber's results used the Littlewood-Richardson rule that mine did not do.

This year I had the urge to return to the triangular matrices and learned of Kostka from MacDonald's book. I decided to implement the procedure for inverting the triangular connection between Frobenius characters and irreducible characters. This process yields the irreducible characters as well as the $a_{k-i}^{k-j}$ matrices. The question is whether or not the monomial approach or the triangular inversion approach is faster. The Young tableaux method was introduced because the monomial approach of Frobenius is "formidable," to quote Hamermesh [6], and the use of Young tableaux is "tractable." Indeed this is why my 1967 paper began with Young tableaux and not with monomials.

## II. Frobenius formula and Kostka matrices.

The remarkable monomial formula of Frobenius connects the Frobenius compound characters with the irreducible characters [6]

(5)

$$\left(\sum_{(\lambda)} \left(\phi_{(l)}^{(\lambda)}\right) \sum_{perm} x_1^{\lambda_1} \ldots x_n^{\lambda_n}\right) \prod_{i<j}^{n}(x_i - x_j)$$
$$= \sum_{(\lambda)} \left(\chi_{(l)}^{(\lambda)}\right) \sum_{P \epsilon S_n} \delta_P P\left(x_1^{\lambda_1+n-1} x_2^{\lambda_2+n-2} \ldots x_n^{\lambda_n}\right)$$

On the right-hand side $P$ denotes a permutation from $S_n$, and $\delta_P$ is $\pm 1$ depending on whether the permutation is even or odd. The sum over $perm$ on the left-hand-side does not mean all permutations but only those producing a distinct monomial; repetitions are not counted (these are sometimes called *simple symmetric forms*). The sums are over character names and the free index, $(l)$ is from the class names. Note that every term on each side of the equation is a monomial of order $n(n + 1)/2$. On the left-hand-side the Vandermonde determinant appears

(6)

$$\prod_{i<j}^{n}(x_i - x_j) = \sum_{P \epsilon S_n} \delta_P P(x_1^{n-1} x_2^{n-2} \ldots x_n^0)$$

In the present notation, the triangularity property can be expressed by

(7)

$$\phi_{(l)}^{(\lambda)} = \sum_{(\mu)} a_{(\mu)}^{(\lambda)} \chi_{(l)}^{(\mu)}$$

$$a_{(\mu)}^{(\lambda)} = \frac{1}{n!} \sum_{(\sigma)} m_{(\sigma)} \phi_{(\sigma)}^{(\lambda)} \chi_{(\sigma)}^{(\mu)}$$

$$a_{(\mu)}^{(\lambda)} = 0 \text{ if } (\mu) < (\lambda)$$

in which the linear ordering of cycle structures was explained in the Introduction. Multiplying the left and right hand sides of eq.(5) by

(8)

$$\frac{1}{n!} \sum_{(l)} m_{(l)} \chi_{(l)}^{(\mu)}$$

and doing the summation yields

(9)

$$\left( \sum_{(\lambda)} \left( a_{(\mu)}^{(\lambda)} \right) \sum_{perm} x_1^{\lambda_1} \ldots x_n^{\lambda_n} \right) \prod_{i<j}^{n} (x_i - x_j)$$
$$= \sum_{P \epsilon S_n} \delta_P P\left( x_1^{\mu_1+n-1} x_2^{\mu_2+n-2} \ldots x_n^{\mu_n} \right)$$

because of eq.(2) and eq.(3). Dividing both sides by the Vandermonde determinant we get an instance of Section 331 from Muir's book on determinants [7].

(10)

$$\sum_{(\lambda)} \left(a_{(\mu)}^{(\lambda)}\right) \sum_{perm} x_1^{\lambda_1} \dots x_n^{\lambda_n} = \frac{\sum_{P \in S_n} \delta_P P\left(x_1^{\mu_1+n-1} x_2^{\mu_2+n-2} \dots x_n^{\mu_n}\right)}{\prod_{i<j}^n (x_i - x_j)}$$

Denote the quotient on the right-hand-side of eq.(10) by Q. Using Section 339 of Muir [7], Q may be rewritten as a determinant of *completely symmetric forms*:

(11)

$$Q = (-1)^{\frac{n(n-1)}{2}} \times$$

$$\text{Det} \begin{bmatrix} (x_1 \dots x_n)^{\mu_1+n-1} & (x_1 \dots x_n)^{\mu_2+n-2} & a & a & a & (x_1 \dots x_n)^{\mu_n} \\ (x_1 \dots x_n)^{\mu_1+n-2} & (x_1 \dots x_n)^{\mu_2+n-3} & a & a & a & (x_1 \dots x_n)^{\mu_n-1} \\ a & a & & a & a & a & a \\ & & & a & a & & a & a \\ (x_1 \dots x_n)^{\mu_1+1} & (x_1 \dots x_n)^{\mu_2} & a & a & a & (x_1 \dots x_n)^{\mu_n-n+2} \\ (x_1 \dots x_n)^{\mu_1} & (x_1 \dots x_n)^{\mu_2-1} & a & a & a & (x_1 \dots x_n)^{\mu_n-n+1} \end{bmatrix}$$

with the notational conventions

(12)

$$(x_1 \dots x_n)^p \equiv 0, \quad p < 0$$
$$(x_1 \dots x_n)^p \equiv 1, \quad p = 0$$
$$(x_1 \dots x_n)^p \equiv \text{the } completely \ symmetric \ form \text{ for } p > 0$$

Example: $(x_1 x_2 x_3)^4 = x_1^4 + x_2^4 + x_3^4 + x_1^3 x_2 + x_1^3 x_3 + x_2^3 x_1 + x_2^3 x_3 + x_3^3 x_1 + x_3^3 x_2 + x_1^2 x_2^2 + x_1^2 x_3^2 + x_2^2 x_3^2 + x_1^2 x_2 x_3 + x_2^2 x_1 x_3 + x_3^2 x_1 x_2$

In this example we see why the short-hand notation of eq.(12) used in eq.(11) is necessary. To make the point emphatically let us compute the values of $a_{(\mu)}^{(\lambda)}$ for $(\mu) = (1^3)$. Eq.(10) and eq.(11) become

(13)

$$\sum_{(\lambda)} a_{(1^3)}^{(\lambda)} \sum_{perm} x_1^{\lambda_1} x_2^{\lambda_2} x_3^{\lambda_3} = -\det \begin{bmatrix} (x_1 x_2 x_3)^3 & (x_1 x_2 x_3)^2 & (x_1 x_2 x_3)^1 \\ (x_1 x_2 x_3)^2 & (x_1 x_2 x_3)^1 & 1 \\ (x_1 x_2 x_3)^1 & 1 & 0 \end{bmatrix}$$

$$= det \begin{bmatrix} (x_1x_2x_3)^3 & (x_1x_2x_3)^2 \\ (x_1x_2x_3)^1 & 1 \end{bmatrix} - (x_1x_2x_3)^1 det \begin{bmatrix} (x_1x_2x_3)^2 & (x_1x_2x_3)^1 \\ (x_1x_2x_3)^1 & 1 \end{bmatrix}$$

$$= (x_1x_2x_3)^3 - 2(x_1x_2x_3)^1(x_1x_2x_3)^2 + ((x_1x_2x_3)^1)^3.$$

Note carefully that there are two types of superscripts, monomial labels for completely symmetric forms as in the example to eq.(12), and ordinary powers as in the 3 in the final term. Just to be clear the last line in eq.(13) becomes

(14)

$$\{x_1^3 + x_2^3 + x_3^3 + x_1^2x_2 + x_1^2x_3 + x_2^2x_1 + x_2^2x_3 + x_3^2x_1 + x_3^2x_2 + x_1x_2x_3\}$$
$$-2\{x_1 + x_2 + x_3\}\{x_1^2 + x_2^2 + x_3^2 + x_1x_2 + x_1x_3 + x_2x_3\}$$
$$+\{x_1 + x_2 + x_3\}\{x_1 + x_2 + x_3\}\{x_1 + x_2 + x_3\}$$

$$= x_1x_2x_3$$

Nearly everything has canceled out leaving just one term. Comparing this with the left-hand-side of eq.(13) we conclude that

(15)

$$a_{(1^3)}^{(1^3)} = 1, a_{(1^3)}^{(2,1)} = 0, a_{(1^3)}^{(3)} = 0$$

This is consistent with eq.(7). Continuing in this manner for $(\mu) = (2,1)$ we get

(16)

$$\sum_{(\lambda)} a_{(2,1)}^{(\lambda)} \sum_{perm} x_1^{\lambda_1} x_2^{\lambda_2} x_3^{\lambda_3} = -det \begin{bmatrix} (x_1x_2x_3)^4 & (x_1x_2x_3)^2 & 1 \\ (x_1x_2x_3)^3 & (x_1x_2x_3)^1 & 0 \\ (x_1x_2x_3)^2 & 1 & 0 \end{bmatrix}$$

$$= -det \begin{bmatrix} (x_1x_2x_3)^3 & (x_1x_2x_3)^1 \\ (x_1x_2x_3)^2 & 1 \end{bmatrix} = -(x_1x_2x_3)^3 + (x_1x_2x_3)^2(x_1x_2x_3)^1$$

$$= x_1^2 x_2 + x_1^2 x_3 + x_2^2 x_1 + x_2^2 x_3 + x_3^2 x_1 + x_3^2 x_2 + 2 x_1 x_2 x_3$$

Comparing this with the left-hand-side of eq.(16) we conclude that

(17)

$$a^{(1^3)}_{(2,1)} = 2, \quad a^{(2,1)}_{(2,1)} = 1, \quad a^{(3)}_{(2,1)} = 0$$

Finally, $(\mu) = (3)$ gives

(18)

$$\sum_{(\lambda)} a^{(\lambda)}_{(3)} \sum_{perm} x_1^{\lambda_1} x_2^{\lambda_2} x_3^{\lambda_3} = -det \begin{bmatrix} (x_1 x_2 x_3)^5 & (x_1 x_2 x_3)^1 & 1 \\ (x_1 x_2 x_3)^4 & 1 & 0 \\ (x_1 x_2 x_3)^3 & 0 & 0 \end{bmatrix}$$

$$= (x_1 x_2 x_3)^3$$
$$= x_1^3 + x_2^3 + x_3^3 + x_1^2 x_2 + x_1^2 x_3 + x_2^2 x_1 + x_2^2 x_3 + x_3^2 x_1 + x_3^2 x_2 + x_1 x_2 x_3$$

From this we conclude that

(19)

$$a^{(1^3)}_{(3)} = 1, \quad a^{(2,1)}_{(3)} = 1, \quad a^{(3)}_{(3)} = 1$$

Thus we have found the Kostka table (matrix) for $S_3$. It may be written in the form of eq.(1)

(20)

| $a^{(\lambda)}_{(\mu)}$ | (3) | (2,1) | $(1^3)$ |
|---|---|---|---|
| (3) | 1 | 0 | 0 |
| (2,1) | 1 | 1 | 0 |
| $(1^3)$ | 1 | 2 | 1 |

The upper matrix index labels the rows and the lower matrix index labels the columns.

Of course it would be preferable to obtain the Kostka matrix inverse directly since we wish to compute the irreducible characters from the Frobenius characters that are known by formula. Return to eq.(10) and eq.(11). Multiply both sides of eq.(10), having already incorporated eq.(11), and sum over the common index to get

(21)

$$\sum_{(\mu)} (a^{-1})^{(\mu)}_{(\nu)} \sum_{(\lambda)} \left(a^{(\lambda)}_{(\mu)}\right) \sum_{perm} x_1^{\lambda_1} \ldots x_n^{\lambda_n} = \sum_{(\mu)} (a^{-1})^{(\mu)}_{(\nu)} Q((\mu))$$

or the wonderful identity

(22)

$$\sum_{perm} x_1^{\nu_1} \ldots x_n^{\nu_n} = \sum_{(\mu)} (a^{-1})^{(\mu)}_{(\nu)} Q((\mu))$$

The $Q((\mu))$ are the determinants labeled by class. In fact we have already computed them above. The right-hand-side of eq.(13) which ended up equal to the right-hand-side of eq.(14) is $Q((1^3))$. The right-hand-side of eq.(16) is $Q((2,1))$ and the right-hand-side of eq.(18) is $Q((3))$. Choose $(\nu) = (3)$. Then eq.(22) reads

(23)

$$x_1^3 + x_2^3 + x_3^3 = (a^{-1})^{(1^3)}_{(3)} Q((1^3)) + (a^{-1})^{(2,1)}_{(3)} Q((2,1)) + (a^{-1})^{(3)}_{(3)} Q((3))$$

which implies

(24)

$$(a^{-1})^{(1^3)}_{(3)} = 1, (a^{-1})^{(2,1)}_{(3)} = -1, (a^{-1})^{(3)}_{(3)} = 1$$

Similarly, for $(\nu) = (2,1)$ we get

(25)

$$x_1^2 x_2 + x_1^2 x_3 + x_2^2 x_1 + x_2^2 x_3 + x_3^2 x_1 + x_3^2 x_2$$
$$= (a^{-1})_{(2,1)}^{(1^3)} Q((1^3)) + (a^{-1})_{(2,1)}^{(2,1)} Q((2,1)) + (a^{-1})_{(2,1)}^{(3)} Q((3))$$

$$(a^{-1})_{(2,1)}^{(1^3)} = -2, (a^{-1})_{(2,1)}^{(2,1)} = 1, (a^{-1})_{(2,1)}^{(3)} = 0$$

and for $(v) = (1^3)$

(26)

$$x_1 x_2 x_3 = (a^{-1})_{(1^3)}^{(1^3)} Q((1^3)) + (a^{-1})_{(1^3)}^{(2,1)} Q((2,1)) + (a^{-1})_{(1^3)}^{(3)} Q((3))$$

$$(a^{-1})_{(1^3)}^{(1^3)} = 1, (a^{-1})_{(1^3)}^{(2,1)} = 0, (a^{-1})_{(1^3)}^{(3)} = 0$$

Thus the inverse Kostka table (matrix) is

(27)

| $(a^{-1})_{(\mu)}^{(\lambda)}$ | (3) | (2,1) | $(1^3)$ |
|---|---|---|---|
| (3) | 1 | 0 | 0 |
| (2,1) | -1 | 1 | 0 |
| $(1^3)$ | 1 | -2 | 1 |

    It is clear that the amount of work required computing the Kostka matrix is the same as required computing its inverse using the monomial method. Since we want the inverse Kostka matrix in order to get the irreducible characters in terms of the Frobenius characters, it is one step easier, a matrix inversion step, to go directly for the inverse using eq.(22). It is also worth noting that the inverse Kostka matrix is obtained by itself without needing to have the Frobenius characters first. The triangular matrix inversion process described in the Introduction and to be executed below requires the Frobenius characters and yields both the Kostka matrix and the irreducible characters all at once. No monomials are involved and the labor is much less as the order of the symmetric group increases.

### III. Using the triangular matrix properties for inversion is faster than using monomials.

The monomial approach quickly leads to enormous expressions as the size of the symmetric group increases. In the example of $S_3$ given in section II, the largest single factor is a completely symmetric form in $Q((3))$ that contains 10 distinct summands. This factor is $(x_1 x_2 x_3)^3$. As the order of the group increases so does the number of and the size of the determinants, which means products of greater numbers of factors. For $S_4$ a corresponding factor is $(x_1 x_2 x_3 x_4)^4$ which when expanded explicitly contains 35 distinct summands. $(x_1 x_2 x_3 x_4 x_5)^5$ contains 126 summands. The combinatorics explodes with increasing $n$. In general, $(x_1 x_2 \ldots x_n)^r$ has $\binom{n+r-1}{r}$ summands. For $n = 8, r = 5$ we get 792 summands and for $n = 12, r = 6$ there are 12,376. A term this large would appear in a 12×12 determinant because $n = 12$ means the group is $S_{12}$ but the number of classes is 77. 77 different determinants would be required to get all the results. Clearly no one is going to do this calculation by hand.

Consider the inversion problem expressed in eq.(1) for $S_3$. It takes the form

(28)

$$\phi^{(3)} = a^{(3)}_{(3)} \chi^{(3)}$$
$$\phi^{(2,1)} = a^{(2,1)}_{(3)} \chi^{(3)} + a^{(2,1)}_{(2,1)} \chi^{(2,1)}$$
$$\phi^{(1^3)} = a^{(1^3)}_{(3)} \chi^{(3)} + a^{(1^3)}_{(2,1)} \chi^{(2,1)} + a^{(1^3)}_{(1^3)} \chi^{(1^3)}$$

The Frobenius compound characters for $S_3$ are given by

(29)

| $m_{(l)}$ | 2 | 3 | 1 |
|---|---|---|---|
| $\phi^{(\lambda)}_{(l)}$ | (3) | (2,1) | $(1^3)$ |
| (3) | 1 | 1 | 1 |
| (2,1) | 0 | 1 | 3 |
| $(1^3)$ | 0 | 0 | 6 |

The rows are named after the character names and the columns are named after the class names. Above a class name is the size of the class. The inversion proceeds from the top row of eq.(28) and uses eq.(3)

(30)

Row 1
$$\chi^{(3)} = \frac{1}{a_{(3)}^{(3)}} \phi^{(3)},$$
$$a_{(3)}^{(3)} = \frac{1}{6}\sum_l m_l \phi^{(3)}(C_l) \chi^{(3)}(C_l) = \frac{1}{a_{(3)}^{(3)}} \frac{1}{6} \sum_l m_l \phi^{(3)}(C_l) \phi^{(3)}(C_l) = \frac{1}{a_{(3)}^{(3)}}$$
Therefore, $a_{(3)}^{(3)} = 1$ and $\chi^{(3)} = \phi^{(3)}$.

Row 2
$$a_{(3)}^{(2,1)} = \frac{1}{6} \sum_l m_l \phi^{(2,1)}(C_l) \chi^{(3)}(C_l) = \frac{1}{6} \sum_l m_l \phi^{(2,1)}(C_l) \phi^{(3)}(C_l) = 1$$
$$\chi^{(2,1)} = \frac{1}{a_{(2,1)}^{(2,1)}} \left( \phi^{(2,1)} - \phi^{(3)} \right)$$
$$a_{(2,1)}^{(2,1)} = \frac{1}{6} \sum_l m_l \phi^{(2,1)}(C_l) \chi^{(2,1)}(C_l)$$
$$= \frac{1}{a_{(2,1)}^{(2,1)}} \frac{1}{6} \sum_l m_l \phi^{(2,1)}(C_l) \left( \phi^{(2,1)}(C_l) - \phi^{(3)}(C_l) \right) = \frac{1}{a_{(2,1)}^{(2,1)}} (2-1)$$
Therefore $a_{(3)}^{(2,1)} = 1$, $a_{(2,1)}^{(2,1)} = 1$, $\chi^{(2,1)} = \left( \phi^{(2,1)} - \phi^{(3)} \right)$

Row 3
$$a_{(3)}^{(1^3)} = \frac{1}{6} \sum_l m_l \phi^{(1^3)}(C_l) \chi^{(3)}(C_l) = \frac{1}{6} \sum_l m_l \phi^{(1^3)}(C_l) \phi^{(3)}(C_l) = 1$$
$$a_{(2,1)}^{(1^3)} = \frac{1}{6} \sum_l m_l \phi^{(1^3)}(C_l) \chi^{(2,1)}(C_l) =$$
$$\frac{1}{6} \sum_l m_l \phi^{(1^3)}(C_l) \left( \phi^{(2,1)}(C_l) - \phi^{(3)}(C_l) \right) = 3 - 1 = 2$$
$$\chi^{(1^3)} = \frac{1}{a_{(1^3)}^{(1^3)}} \left( \phi^{(1^3)} - \phi^{(3)} - 2\left( \phi^{(2,1)} - \phi^{(3)} \right) \right)$$
$$= \frac{1}{a_{(1^3)}^{(1^3)}} \left( \phi^{(1^3)} + \phi^{(3)} - 2\phi^{(2,1)} \right)$$

$$a^{(1^3)}_{(1^3)} = \frac{1}{6}\sum_l m_l \phi^{(1^3)}(C_l)\chi^{(1^3)}(C_l) = \frac{1}{a^{(1^3)}_{(1^3)}}\frac{1}{6}\sum_l m_l \phi^{(1^3)}(C_l) \times$$
$$\left(\phi^{(1^3)}(C_l) + \phi^{(3)}(C_l) - 2\phi^{(2,1)}(C_l)\right) =$$
$$\frac{1}{a^{(1^3)}_{(1^3)}}(6 + 1 - 6)$$

Therefore $a^{(1^3)}_{(3)} = 1, a^{(1^3)}_{(2,1)} = 2, a^{(1^3)}_{(1^3)} = 1, \chi^{(1^3)} = \phi^{(1^3)} + \phi^{(3)} - 2\phi^{(2,1)}$

From these results we can construct the irreducible character table.

(31)

| $m_{(l)}$ | 2 | 3 | 1 |
|---|---|---|---|
| $\chi^{(\lambda)}_{(l)}$ | (3) | (2,1) | $(1^3)$ |
| (3) | 1 | 1 | 1 |
| (2,1) | -1 | 0 | 2 |
| $(1^3)$ | 1 | -1 | 1 |

The amount of work done here compares favorably with the monomial calculation in which 3 3×3 determinants were required. Now consider what must be done using both methods if we work out $S_5$. Eq.(22) requires 7 Q-determinants each of size 5×5. The elements in these determinants are completely symmetric monomials of the form $(x_1 x_2 x_3 x_4 x_5)^p$ in which, according to eq.(11) the largest value of $p$ is $\mu_1 + 5 - 1$ in the upper left-hand corner of the determinant and the smallest value (most negative) of $p$ is $\mu_n - 5 + 1$ in the lower right-hand corner, for class $(\mu)$. When a determinant of size 5×5 is evaluated there are 5! (120) products of 5 factors each that sum to the determinant. Many of the lower right-hand corner elements are zero because $p < 0$ or 1 because $p = 0$. That helps but not that much for larger $n$. Look back at eq.(13), eq.(16) and eq.(18) to see this for $S_3$. There, 3 3×3 determinants had to be computed but because of the ones and zeros only 4 2×2 determinants were non-zero. Instead of a possible 18 products only 8 were needed, and one of these was zero while four others were multiplications by one. Nevertheless as $n$ increases the number of multiplications grows rapidly and the multiplicands get much longer. Going from 3 to 5 may make the point manifest.

# III. Solving $S_5$ by the triangle method.

To begin solving $S_5$ we need the Frobenius compound characters for both the monomial solution and the triangular solution. A short digression to compute these compound characters is presented here.

We follow Hamermesh [6] and Littlewood [8] (used also by Hamermesh). Their presentations illustrate how confusion is created by notational choices. I will use a notation that I believe eliminates the confusions. Let $G$ be the symmetric group $S_n$. The order of $G$ is $g = n!$. A sub-group of $G$ can be constructed from any cycle structure of $G$ given by $(\lambda_1, \lambda_2, \ldots, \lambda_n)$ in which $\lambda_1 + \lambda_2 + \cdots + \lambda_n = n$ and $\lambda_1 \geq \lambda_2 \geq \cdots \geq \lambda_n$. Take the $n$ integers $1, 2, \ldots, n$ and place any $\lambda_1$ of them in set 1, any $\lambda_2$ of the rest of them in set 2, any $\lambda_3$ of the rest of them in set 3, etc. until done. The order of the integers in these sets is irrelevant so long as no two sets have any integers in common. Take set 1 and construct the symmetric group on the $\lambda_1$ integers forming the group $S_{\lambda_1}$ of order $\lambda_1!$. Do the same for each set. The direct product group, $H_{(\lambda)} = S_{\lambda_1} \times S_{\lambda_2} \times \ldots \times S_{\lambda_n}$ is a sub-group of $S_n$. It has order $\lambda_1! \lambda_2! \ldots \lambda_n!$.

Consider the class $(l)$ of $S_n$. The cycle structure of $(l)$ is written $\left(l_1^{m_1}, l_2^{m_2}, l_3^{m_3}, \ldots, l_n^{m_n}\right)$ which means there are $m_1$ cycles of length $l_1 = 1$, $m_2$ cycles of length $l_2 = 2$, $m_3$ cycles of length $l_3 = 3$ etc., such that $m_1 l_1 + m_2 l_2 + \cdots + m_n l_n = n$. The number of permutations in $S_n$ from class $(l)$, denoted by $g_{(l)}$, is

(32)

$$g_{(l)} = \frac{n!}{l_1^{m_1} m_1! \, l_2^{m_2} m_2! \ldots l_n^{m_n} m_n!}$$

The quantity $h_{(l)}$ is the number of permutations in $H_{(\lambda)}$ from class $(l)$ of $S_n$. In order for a permutation in $H_{(\lambda)}$ to contribute to $h_{(l)}$ it is necessary that it contain
$m_1 \, l_1 - cycles, m_2 \, l_2 - cycles, m_3 \, l_3 - cycles, \ldots, m_n \, l_n - cycles$
Since $H_{(\lambda)}$ is a direct product of $S_{\lambda_i}$'s, $i = 1, 2, \ldots, n$, a permutation with factors in $S_{\lambda_1}$ must have
$m_{1,1} \, l_1 - cycles, m_{2,1} \, l_2 - cycles, m_{3,1} \, l_3 - cycles, \ldots, m_{n,1} \, l_n - cycles$

and with factors in $S_{\lambda_2}$ must have

$m_{1,2}\ l_1-cycles, m_{2,2}\ l_2-cycles, m_{3,2}\ l_3-cycles,\ldots, m_{n,2}\ l_n-cycles$

and so forth, and with factors in $S_{\lambda_n}$ must have

$m_{1,n}\ l_1-cycles, m_{2,n}\ l_2-cycles, m_{3,n}\ l_3-cycles,\ldots, m_{n,n}\ l_n-cycles$

provided

(33)

$$\sum_{i=1}^{n} m_{j,i} = m_j$$

$$\sum_{j=1}^{n} m_{j,i}\ l_j = \lambda_i$$

The first identity holds for each $j = 1,2,\ldots,n$ and the second identity holds for each $i = 1,2,\ldots,n$. Note that $l_j \equiv j$. Thus, if any $m_j = 0$ then all $m_{j,i} = 0$ for the same $j$, or if any $\lambda_i = 0$ then all $m_{j,i} = 0$ for the same $i$. The number of permutations in $S_{\lambda_i}$ having cycle structure $(l_1^{m_{1,i}}, l_2^{m_{2,i}}, l_3^{m_{3,i}}, \ldots, l_n^{m_{n,i}})$ is

$$\frac{\lambda_i!}{l_1^{m_{1,i}} m_{1,i}!\ l_2^{m_{2,i}} m_{2,i}! \ldots l_n^{m_{n,i}} m_{n,i}!}$$

Therefore, the quantity $h_{(l)}$ is given by

(34)

$$h_{(l)} = \sum \prod_i \frac{\lambda_i!}{l_1^{m_{1,i}} m_{1,i}!\ l_2^{m_{2,i}} m_{2,i}! \ldots l_n^{m_{n,i}} m_{n,i}!}$$

in which the summation is over all solutions to eq.(33).

The Frobenius compound characters are defined by

(35)

$$\phi_{(l)}^{(\lambda)} \equiv \frac{gh_{(l)}}{g_{(l)}h} = \sum_j \prod \frac{m_j!}{\prod_i m_{j,i}!}$$

in which the summation is over all solutions to eq.(33). A remarkable amount of simplification by cancellation has occurred. The character symbol on the left-hand side of eq.(35) depends on the character name cycle structure $(\lambda)$ and on the class structures $(l)$. On the right-hand-side we see no sign of either. Both are explicit in the right-hand-sides of eqs.(33).

For $S_5$ the Frobenius character table is

(36)

| $m_{(l)}$ | 24 | 30 | 20 | 20 | 15 | 10 | 1 |
|---|---|---|---|---|---|---|---|
| $\phi_{(l)}^{(\lambda)}$ | $(5)$ | $(4,1)$ | $(3,2)$ | $(3,1^2)$ | $(2^2,1)$ | $(2,1^3)$ | $(1^5)$ |
| $(5)$ | 1 | 1 | 1 | 1 | 1 | 1 | 1 |
| $(4,1)$ | 0 | 1 | 0 | 2 | 1 | 3 | 5 |
| $(3,2)$ | 0 | 0 | 1 | 1 | 2 | 4 | 10 |
| $(3,1^2)$ | 0 | 0 | 0 | 2 | 0 | 6 | 20 |
| $(2^2,1)$ | 0 | 0 | 0 | 0 | 2 | 6 | 30 |
| $(2,1^3)$ | 0 | 0 | 0 | 0 | 0 | 6 | 60 |
| $(1^5)$ | 0 | 0 | 0 | 0 | 0 | 0 | 120 |

The character name $(\lambda)$ also labels the rows and the class label $(l)$ labels the columns. I will demonstrate the applicability of eq.(35) by computing a few of the entries in the table.

$(\lambda) = (3,2): \lambda_1 = 3, \lambda_2 = 2$
$(l) = (2,1^3): m_1 = 3, m_2 = 1$
$m_{1,1} + 2m_{2,1} + 3m_{3,1} + 4m_{4,1} + 5m_{5,1} = 3$
$m_{1,2} + 2m_{2,2} + 3m_{3,2} + 4m_{4,2} + 5m_{5,2} = 2$
$m_{1,1} + m_{1,2} + m_{1,3} + m_{1,4} + m_{1,5} = 3$
$m_{2,1} + m_{2,2} + m_{2,3} + m_{2,4} + m_{2,5} = 1$
$\{m_{1,1} = 3, m_{2,2} = 1\}$ or $\{m_{1,1} = 1, m_{1,2} = 2, m_{2,1} = 1\}$
$\frac{3!1!}{3!1!} = 1$ and $\frac{3!1!}{1!2!1!} = 3$
$3 + 1 = 4$

$(\lambda) = (2^2, 1)$: $\lambda_1 = 2, \lambda_2 = 2, \lambda_3 = 1$
$(l) = (1^5),$ : $m_1 = 5$
$m_{1,1} + 2m_{2,1} + 3m_{3,1} + 4m_{4,1} + 5m_{5,1} = 2$
$m_{1,2} + 2m_{2,2} + 3m_{3,2} + 4m_{4,2} + 5m_{5,2} = 2$
$m_{1,3} + 2m_{2,3} + 3m_{3,3} + 4m_{4,3} + 5m_{5,3} = 1$
$m_{1,1} + m_{1,2} + m_{1,3} + m_{1,4} + m_{1,5} = 5$
$\{m_{1,1} = 2, m_{1,2} = 2, m_{1,3} = 1\}$
$$\frac{5!}{2!\,2!\,1!} = \mathbf{30}$$

$(\lambda) = (3,2)$: $\lambda_1 = 3, \lambda_2 = 2$
$(l) = (4,1)$: $m_1 = 1, m_4 = 1$
$m_{1,1} + 2m_{2,1} + 3m_{3,1} + 4m_{4,1} + 5m_{5,1} = 3$
$m_{1,2} + 2m_{2,2} + 3m_{3,2} + 4m_{4,2} + 5m_{5,2} = 2$
$m_{1,1} + m_{1,2} + m_{1,3} + m_{1,4} + m_{1,5} = 1$
$m_{4,1} + m_{4,2} + m_{4,3} + m_{4,4} + m_{4,5} = 1$
no solutions
**0**

In the last example, $m_2$, $m_3$, and $m_5$ are zero. Therefore any term of the form $m_{j,i}$ with $j = 2, 3, or\ 5$ and $i = 1, 2, 3, 4, 5$ must vanish. What is left for the two $\lambda$-equations has no solutions. Indeed 21 of the entries in the table vanish for similar reasons. With a little practice one can fill out the entire Frobenius character table fairly easily.

The triangular matrix solution process begins with the $S_5$ realization of eq.(1).

(37)

$\phi^7 = a_7^7 \chi^7$
$\phi^6 = a_7^6 \chi^7 + a_6^6 \chi^6$
$\phi^5 = a_7^5 \chi^7 + a_6^5 \chi^6 + a_5^5 \chi^5$
$\phi^4 = a_7^4 \chi^7 + a_6^4 \chi^6 + a_5^4 \chi^5 + a_4^4 \chi^4$
$\phi^3 = a_7^3 \chi^7 + a_6^3 \chi^6 + a_5^3 \chi^5 + a_4^3 \chi^4 + a_3^3 \chi^3$
$\phi^2 = a_7^2 \chi^7 + a_6^2 \chi^6 + a_5^2 \chi^5 + a_4^2 \chi^4 + a_3^2 \chi^3 + a_2^2 \chi^2$
$\phi^1 = a_7^1 \chi^7 + a_6^1 \chi^6 + a_5^1 \chi^5 + a_4^1 \chi^4 + a_3^1 \chi^3 + a_2^1 \chi^2 + a_1^1 \chi^1$

The notational correspondence is 7, 6, 5, 4, 3, 2, 1 with $(5), (4,1), (3,2), (3, 1^2), (2^2, 1), (2, 1^3), (1^5)$. A *bracket* notation, <brac|ket>, is introduced to further simplify the notation so that

$$\frac{1}{n!}\sum_\sigma m_\sigma \phi^j(C_\sigma)\chi^i(C_\sigma) \to <\phi^j|\chi^i>$$

$$\frac{1}{n!}\sum_\sigma m_\sigma \phi^j(C_\sigma)\phi^i(C_\sigma) \to <\phi^j|\phi^i>$$

Row 1
$$\chi^7 = \frac{1}{a_7^7}\phi^7$$

$$a_7^7 = <\phi^7|\chi^7> = \frac{1}{a_7^7}<\phi^7|\phi^7> = \frac{1}{a_7^7}$$

$$a_7^7 = 1, \chi^7 = \phi^7$$

Row 2
$$a_7^6 = <\phi^6|\chi^7> = <\phi^6|\phi^7> = 1$$
$$\chi^6 = \frac{1}{a_6^6}(\phi^6 - \phi^7)$$

$$a_6^6 = <\phi^6|\chi^6> = \frac{1}{a_6^6}<\phi^6|(\phi^6 - \phi^7)> = \frac{1}{a_6^6}(2-1)$$

$$a_7^6 = 1, a_6^6 = 1, \chi^6 = (\phi^6 - \phi^7)$$

Row 3
$$a_7^5 = <\phi^5|\chi^7> = <\phi^5|\phi^7> = 1$$
$$a_6^5 = <\phi^5|\chi^6> = <\phi^5|(\phi^6 - \phi^7)> = (2-1)$$
$$a_5^5 = <\phi^5|\chi^5> = \frac{1}{a_5^5}<\phi^5|(\phi^5 - \phi^6)> = \frac{1}{a_5^5}(3-2)$$

$$a_7^5 = 1, a_6^5 = 1, a_5^5 = 1, \chi^5 = (\phi^5 - \phi^6)$$

Row 4
$$a_7^4 = <\phi^4|\chi^7> = <\phi^4|\phi^7> = 1$$
$$a_6^4 = <\phi^4|\chi^6> = <\phi^4|(\phi^6 - \phi^7)> = (3-1)$$
$$a_5^4 = <\phi^4|\chi^5> = <\phi^4|(\phi^5 - \phi^6)> = (4-3)$$

$$a_4^4 = <\phi^4|\chi^4> = \frac{1}{a_4^4}<\phi^4|(\phi^4 - \phi^7 - 2(\phi^6 - \phi^7) - (\phi^5 - \phi^6))> =$$
$$\frac{1}{a_4^4}<\phi^4|(\phi^4 + \phi^7 - \phi^6 - \phi^5)> = \frac{1}{a_4^4}(7 + 1 - 3 - 4)$$

$a_7^4 = 1, a_6^4 = 2, a_5^4 = 1, a_4^4 = 1,$
$\chi^4 = (\phi^4 + \phi^7 - \phi^6 - \phi^5)$

Row 5
$a_7^3 = <\phi^3|\chi^7> = <\phi^3|\phi^7> = 1$
$a_6^3 = <\phi^3|\chi^6> = <\phi^3|(\phi^6 - \phi^7)> = (3 - 1)$
$a_5^3 = <\phi^3|\chi^5> = <\phi^3|(\phi^5 - \phi^6)> = (5 - 3)$
$a_4^3 = <\phi^3|\chi^4> = <\phi^3|(\phi^4 + \phi^7 - \phi^6 - \phi^5)> = (8 + 1 - 3 - 5)$
$a_3^3 = <\phi^3|\chi^3> = \frac{1}{a_3^3}<\phi^3|(\phi^3 - \phi^7 - 2(\phi^6 - \phi^7))>$

$-\frac{1}{a_3^3}<\phi^3|(2(\phi^5 - \phi^6) + (\phi^4 + \phi^7 - \phi^6 - \phi^5))>$
$= \frac{1}{a_3^3}<\phi^3|(\phi^3 + \phi^6 - \phi^5 - \phi^4)> = \frac{1}{a_3^3}(11 + 3 - 5 - 8)$

$a_7^3 = 1, a_6^3 = 2, a_5^3 = 2, a_4^3 = 1, a_3^3 = 1,$
$\chi^3 = (\phi^3 + \phi^6 - \phi^5 - \phi^4)$

Row 6
$a_7^2 = <\phi^2|\chi^7> = <\phi^2|\phi^7> = 1$
$a_6^2 = <\phi^2|\chi^6> = <\phi^2|(\phi^6 - \phi^7)> = (4 - 1)$
$a_5^2 = <\phi^2|\chi^5> = <\phi^2|(\phi^5 - \phi^6)> = (7 - 4)$
$a_4^2 = <\phi^2|\chi^4> = <\phi^2|(\phi^4 + \phi^7 - \phi^6 - \phi^5)> = (13 + 1 - 4 - 7)$
$a_3^2 = <\phi^2|\chi^3> = <\phi^2|(\phi^3 + \phi^6 - \phi^5 - \phi^4)> = (18 + 4 - 7 - 13)$
$a_2^2 = <\phi^2|\chi^2> = \frac{1}{a_2^2}<\phi^2|(\phi^2 - \phi^7 - 3(\phi^6 - \phi^7))>$

$-\frac{1}{a_2^2}<\phi^2|(3(\phi^5 - \phi^6) + 3(\phi^4 + \phi^7 - \phi^6 - \phi^5))>$
$-\frac{1}{a_2^2}<\phi^2|2(\phi^3 + \phi^6 - \phi^5 - \phi^4)>$

$$= \frac{1}{a_2^2} < \phi^2 | (\phi^2 - \phi^7 + \phi^6 + 2\phi^5 - \phi^4 - 2\phi^3) >$$

$$= \frac{1}{a_2^2} (33 - 1 + 4 + 14 - 13 - 36)$$

$a_7^2 = 1, a_6^2 = 3, a_5^2 = 3, a_4^2 = 3, a_3^2 = 2, a_2^2 = 1,$
$\chi^2 = (\phi^2 - \phi^7 + \phi^6 + 2\phi^5 - \phi^4 - 2\phi^3)$

Row 7
$a_7^1 = < \phi^1 | \chi^7 > = < \phi^1 | \phi^7 > = 1$
$a_6^1 = < \phi^1 | \chi^6 > = < \phi^1 | (\phi^6 - \phi^7) > = (5 - 1)$
$a_5^1 = < \phi^1 | \chi^5 > = < \phi^1 | (\phi^5 - \phi^6) > = (10 - 5)$
$a_4^1 = < \phi^1 | \chi^4 > = < \phi^1 | (\phi^4 + \phi^7 - \phi^6 - \phi^5) > = (20 + 1 - 5 - 10)$
$a_3^1 = < \phi^1 | \chi^3 > = < \phi^1 | (\phi^3 + \phi^6 - \phi^5 - \phi^4) > = (30 + 5 - 10 - 20)$
$a_2^1 = < \phi^1 | \chi^2 > = < \phi^1 | (\phi^2 - \phi^7 + \phi^6 + 2\phi^5 - \phi^4 - 2\phi^3) >$
$$= (60 - 1 + 5 + 20 - 20 - 60)$$
$a_1^1 = < \phi^1 | \chi^1 > = \frac{1}{a_1^1} < \phi^1 | (\phi^1 - \phi^7 - 4(\phi^6 - \phi^7)) >$

$$-\frac{1}{a_1^1} < \phi^1 | (5(\phi^5 - \phi^6) + 6(\phi^4 + \phi^7 - \phi^6 - \phi^5)) >$$

$$-\frac{1}{a_1^1} < \phi^1 | (5(\phi^3 + \phi^6 - \phi^5 - \phi^4)$$
$$+ 4(\phi^2 - \phi^7 + \phi^6 + 2\phi^5 - \phi^4 - 2\phi^3)) >$$
$$= \frac{1}{a_1^1} < \phi^1 | (\phi^1 + \phi^7 - 2\phi^6 - 2\phi^5 + 3\phi^4 + 3\phi^3 - 4\phi^2)$$
$$= \frac{1}{a_1^1} (120 + 1 - 10 - 20 + 60 + 90 - 240)$$

$a_7^1 = 1, a_6^1 = 4, a_5^1 = 5, a_4^1 = 6, a_3^1 = 5, a_2^1 = 4, a_1^1 = 1,$
$\chi^1 = (\phi^1 + \phi^7 - 2\phi^6 - 2\phi^5 + 3\phi^4 + 3\phi^3 - 4\phi^2)$

### IV. A nice surprise

The last two lines of each Row calculation contain a catalog of the results for that Row. The $a_i^j$'s are the Kostka matrix elements and the $\chi^j$'s are the irreducible characters indirectly given by linear expansions in terms of the Frobenius characters. The coefficients of these expansions are precisely the matrix elements of the inverse Kostka matrix! We get them for free. Reading from the Row results we get:

(38)

Kostka matrix (table):

| 1 | 0 | 0 | 0 | 0 | 0 | 0 |
|---|---|---|---|---|---|---|
| 1 | 1 | 0 | 0 | 0 | 0 | 0 |
| 1 | 1 | 1 | 0 | 0 | 0 | 0 |
| 1 | 2 | 1 | 1 | 0 | 0 | 0 |
| 1 | 2 | 2 | 1 | 1 | 0 | 0 |
| 1 | 3 | 3 | 3 | 2 | 1 | 0 |
| 1 | 4 | 5 | 6 | 5 | 4 | 1 |

Inverse Kostka matrix (table):

| 1 | 0 | 0 | 0 | 0 | 0 | 0 |
|---|---|---|---|---|---|---|
| -1 | 1 | 0 | 0 | 0 | 0 | 0 |
| 0 | -1 | 1 | 0 | 0 | 0 | 0 |
| 1 | -1 | -1 | 1 | 0 | 0 | 0 |
| 0 | 1 | -1 | -1 | 1 | 0 | 0 |
| -1 | 1 | 2 | -1 | -2 | 1 | 0 |
| 1 | -2 | -2 | 3 | 3 | -4 | 1 |

Using the Frobenius character table in eq.(36) and the inverse Kostka matrix we get the irreducible character table for $S_5$.

| $m_{(l)}$ | 24 | 30 | 20 | 20 | 15 | 10 | 1 |
|---|---|---|---|---|---|---|---|
| $\chi^{(\lambda)}_{(l)}$ | (5) | (4,1) | (3,2) | $(3,1^2)$ | $(2^2,1)$ | $(2,1^3)$ | $(1^5)$ |
| (5) | 1 | 1 | 1 | 1 | 1 | 1 | 1 |
| (4,1) | -1 | 0 | -1 | 1 | 0 | 2 | 4 |
| (3,2) | 0 | -1 | 1 | -1 | 1 | 1 | 5 |
| $(3,1^2)$ | 1 | 0 | 0 | 0 | -2 | 0 | 6 |
| $(2^2,1)$ | 0 | 1 | -1 | -1 | 1 | -1 | 5 |
| $(2,1^3)$ | -1 | 0 | 1 | 1 | 0 | -2 | 4 |
| $(1^5)$ | 1 | -1 | -1 | 1 | 1 | -1 | 1 |

Also notice that the Kerber identities mentioned near the end of section I have been justified in these constructions.

## V. Solving $S_5$ by the monomial method

The construction of irreducible characters for $S_5$ using the monomial method is most beautifully rendered by eq.(22). All the work is contained in computing the 7 determinants, $Q((\mu))$, on the right-hand-side, one for each of 7 cycle structures. In this case these determinants are for 5×5 matrices the elements of which are completely symmetric forms on 5 symbols. The generic structure for the 7 cycle structures produced by 5 variables is:

(39)

$$\begin{vmatrix} (x_1 \ldots x_5)^{\mu_1+4} & (x_1 \ldots x_5)^{\mu_2+3} & (x_1 \ldots x_5)^{\mu_3+2} & (x_1 \ldots x_5)^{\mu_4+1} & (x_1 \ldots x_5)^{\mu_5-0} \\ (x_1 \ldots x_5)^{\mu_1+3} & (x_1 \ldots x_5)^{\mu_2+2} & (x_1 \ldots x_5)^{\mu_3+1} & (x_1 \ldots x_5)^{\mu_4+0} & (x_1 \ldots x_5)^{\mu_5-1} \\ (x_1 \ldots x_5)^{\mu_1+2} & (x_1 \ldots x_5)^{\mu_2+1} & (x_1 \ldots x_5)^{\mu_3} & (x_1 \ldots x_5)^{\mu_4-1} & (x_1 \ldots x_5)^{\mu_5-2} \\ (x_1 \ldots x_5)^{\mu_1+1} & (x_1 \ldots x_5)^{\mu_2+0} & (x_1 \ldots x_5)^{\mu_3-1} & (x_1 \ldots x_5)^{\mu_4-2} & (x_1 \ldots x_5)^{\mu_5-3} \\ (x_1 \ldots x_5)^{\mu_1+0} & (x_1 \ldots x_5)^{\mu_2-1} & (x_1 \ldots x_5)^{\mu_3-2} & (x_1 \ldots x_5)^{\mu_4-3} & (x_1 \ldots x_5)^{\mu_5-4} \end{vmatrix}$$

If $(\mu) = (3, 1^2) = (3,1,1,0,0)$ then $\mu_2 - 1 = 0$, $\mu_3 - 1 = 0$, $\mu_4 = 0$ and $\mu_5 = 0$. The rules in eq.(12) convert the generic matrix in eq.(39) into the particular case

(40)

$$\begin{vmatrix} (x_1 \ldots x_5)^{\mu_1+4} & (x_1 \ldots x_5)^{\mu_2+3} & (x_1 \ldots x_5)^{\mu_3+2} & (x_1 \ldots x_5)^{\mu_4+1} & 1 \\ (x_1 \ldots x_5)^{\mu_1+3} & (x_1 \ldots x_5)^{\mu_2+2} & (x_1 \ldots x_5)^{\mu_3+1} & 1 & \\ (x_1 \ldots x_5)^{\mu_1+2} & (x_1 \ldots x_5)^{\mu_2+1} & (x_1 \ldots x_5)^{\mu_3} & a & 0 \\ (x_1 \ldots x_5)^{\mu_1+1} & (x_1 \ldots x_5)^{\mu_2+0} & 1 & a & 0 \\ (x_1 \ldots x_5)^{\mu_1+0} & 1 & 0 & a & 0 \end{vmatrix}$$

All of the ones are evident and a few of the zeros. The remaining empty spaces are just zeros. Thus the determinant that would include 120 entries if they were all non-zero is only as hard to compute as would be a 3×3 determinant. Even that one resulting above is really 2 2×2 determinants. Other examples are more or less complicated depending on the cycle structure. Unlike the case earlier for $S_3$, where 3 variable monomials were involved, here we have 5 variable monomials and the biggest completely symmetric form in the reduced 3×3 determinant above is $(x_1 \ldots x_5)^{\mu_1+2} = (x_1 \ldots x_5)^5$. As was shown earlier the number of summands in this one case is

(41)

$$\binom{9}{5} = 126$$

That is big number of summands but in the matrix above it gets multiplied by two ones. A worse term is the product $(x_1 \ldots x_5)^{\mu_1+0} \times (x_1 \ldots x_5)^{\mu_2+0} \times (x_1 \ldots x_5)^{\mu_3} = (x_1 \ldots x_5)^3 \times (x_1 \ldots x_5)^1 \times (x_1 \ldots x_5)^1$ which has

(42)

$$\binom{7}{3} \times \binom{5}{1} \times \binom{5}{1} = 35 \times 5 \times 5 = 875$$

That's a lot of terms. For $S_3$ the biggest single monomial had just 10 summands.

It should now be clear that the monomial approach involves a great deal of labor for small $n$ even if the determinants are guaranteed to have a lower right-hand corner full of zeros and ones. These sorts of computations are more suited to computer analysis than to pen and paper. For this reason many other researchers have developed algebraic techniques that do not involve explicit display of tables as was done here. My motivation was to understand better the triangular relationship that I discovered many years ago using only algebra in a study of idempotents. I now know that I was looking at Carl Kostka's matrices. I know two different ways to construct them, the monomial way in eq.(22) and the triangular way in Section III. The triangular way wins a speed race. As for finally being able to give my triangular matrices the name Kostka, I recall that Juliet mused: "a rose by any other name would smell as sweet." Richard Feynman noted the difference between knowing the name of something and knowing something.

## VI. Conclusion

In [1] I asked for a combinatorial interpretation for what I now know are Kostka matrices. Such an interpretation is given by eq.(22). The left-hand-side is expressed solely in terms of simple symmetric forms while the right-hand-side is made up of completely symmetric forms. The inverse Kostka matrix provides the coupling. In Section II this was explicitly exhibited for $S_3$. The example of a completely symmetric form is given below eq.(12). What is shown is the completely symmetric form for three variables and a total power of 4 that yields terms that may be labeled (4), (3,1), (2,2), or (2,1,1). The simply symmetric forms correspond to the terms with just one of these power-class labels. Thus the example can be broken into 4 simply symmetric monomials. In [1] the Kostka matrix provides the interconversion of Frobenius compound characters and irreducible characters. It's interpretation is given by eq.(22) in which it provides the interconversion of simply symmetric forms and completely symmetric forms. Moreover there is no sign of the symmetric group characters in this monomial expression. Perhaps this is the fundamental interpretation underlying the Kostka matrices.

## Apologia

To those readers who know all of this already I request your forbearance. My personal goal was to finish something I started 51 years ago. I made several discoveries even though others may have made them earlier. It was still a pleasure to do so. I thank Carl Kostka for patiently waiting 136 years for me, justifiably, to put his name on my matrices. In the future I will replace my $a_i^j$'s with $K_\nu^\mu$'s.

[4] H. O. Foulkes, "Recurrences for the Characters of the Symmetric Group," *Discrete Mathematics* **21** 137-144 (1978).

[5] A. Kerber, "On a Paper of Fox about a Method for Calculating the Ordinary Irreducible Characters of Symmetric Groups," *Journal of Combinatorial Theory*, **6** 90-93 (1969).

[6] M. Hamermesh, *Group Theory*, p. 197 (Addison-Wesley Publishing Co., Reading, Massachusetts, 1962).

[7] T. Muir, *A Treatise on the Theory of Determinants,* (Longmans, Green and Co., London, 1933).

[8] D. E. Littlewood, *The Theory of Group Characters and Matrix Representations of Groups*, Chapter VIII (Oxford University Press, London, 1950)